\documentclass[reqno]{amsart}
\usepackage{hyperref}
\usepackage{graphicx}

\numberwithin{equation}{section}
\newtheorem{theorem}{Theorem}[section]

\newtheorem{definition}[theorem]{Definition}

\newtheorem{remark}[theorem]{Remark}

\begin{document}
\title{Homogenization of a Double Porosity Model in Deformable Media}
\author{Abdelhamid AINOUZ} 
\address{Abdelhamid AINOUZ \newline
Laboratoire AMNEDP \\
Mathematics Department, USTHB\\
BP 32 ElAlia, Bab-Ezzouar, Algiers, Algeria.
}
\email{aainouz@usthb.dz}
\thanks{}

\begin{abstract}
The paper addresses the homogenization of a family of micro-models for the
flow of a slightly compressible fluid in a poroelastic matrix containing
periodically distibuted poroelastic inclusions, with \ low permeabilities
and with imperfect contact on the interface. The micro-models are based on
Biot's system for consolidation processes in each phase, with interfacial
barrier formulation. Using the two-scale convergence technique, it is shown
that the derived system is a general model of that proposed by Aifantis,
plus an extra memory term.
\end{abstract}
\maketitle

\allowdisplaybreaks

\section{Introduction\label{s1}}

The interaction between fluid flow and solid deformation in porous media is
of great importance in petroleum engineering and geomechanics, biosciences,
chemical processes\ and many industrial applications \cite{cou, cow, wang}.

Some type of porous rocks, like aquifers and petroleum reservoir systems,
may contain fractures. It is known that flows in such media occur mainly in
the fracture region and the dominant fluid storage is in the matrix blocks.
In this situation, reservoir possesses two porous structures, one related to
the matrix, and the other related to fractures. This notion of double
porosity/permeability has first been introduced by Barenblatt, Zheltov and
Kochina \cite{bzk} to model the flow of a slightly compressible fluid within
naturally fractured porous media. The proposed model is a system of two
partial differential equations in a two-medium description, with Darcy's law
in each phase, plus exchange extra-terms representing the interfacial
coupling that results from the interaction, at the micro-scale, between the
two phases, see (\ref{06})-(\ref{07}) below. This was derived under the main
assumption that the fluid pressure is uniformly distributed at the surface
of each phase.

Generally, fractured rock formations present at the micro-scale high degrees
of heterogeneity and permeability is mainly determined by the size of the
pores and the connectedness of the fractures system. So any mathematical
modeling of fluid flow in such porous media must take into account the rapid
spatial variation of the phenomenological parameters. Furthermore, from the
numerical point of view, modeling of such systems at the local scale yields
a huge number of discretized equations, so computations will be fastidious
and intractable. To deal with such highly heterogeneous domains, the idea is
to replace the medium by an effective one. Homogenization techniques, like
two-scale convergence method, have been used to rigorously derive an
effective double-porosity model for the Barenblatt, Zheltov and Kochina
(BZK) system, see for instance H. Ene and D. Polisevski \cite{ep}. However,
this model does not take into account the elastic behavior of the solid. In
fact, a rise in pore pressure of the fluid produces a dilation of the solid
mass. On the other hand, compression of the medium will increase pore
pressure. This coupled pressure-deformation was first introduced by Terzaghi
\cite{ter1} in the one-dimensional setting and gave the first soil
consolidation problem for a homogeneous elastic porous medium. Later, M. A.
Biot \cite{biot} has developed in the multidimensional setting a linear
theoretical analysis for the behavior of a fluid saturated poroelastic
medium. The model was based on macroscopic description of the
phenomenological and physical quantities where the representative volume
element is described as the superposition of a particle of fluid and a
particle of solid. Assuming that microstructures are periodically
distributed and that the pore scale is very small compared to the
macroscopic scale, a two-scale asymptotic expansion technique can be used to
rigorously justify this Biot's model. The microscopic models are based on
the linear elasticity equations in the skeleton and on the Stokes equations
in the fluid with appropriate transmission conditions. For more details, we
refer the reader to the earlier work by Auriault and Sanchez-Palencia \cite%
{as}.

Because of the coupling between the deformation and fluid pressure in double
porosity rocks, which must be understood in order to predict reservoir or
aquifer behavior, the concept of double porosity has been developed by E.C.
Aifantis \cite{aif}\ to model oil flow in porous elastic rocks. More
precisely, E. C. Aifantis gave a phenomenological model for flow of a weakly
compressible fluid in a complex and heterogeneous medium where a system of
partial differential equations is given and generalizing Biot's
consolidation model by taking into account the basic physics of flow through
fractured media with interscale couplings. The proposed model reads as
follows:
\begin{eqnarray}
-\mu \Delta \mathbf{u}-(\lambda +\mu )\nabla (\mathrm{div}\mathbf{u})+\alpha
_{1}\nabla p_{1}+\alpha _{2}\nabla p_{2} &=&\mathbf{f}\text{, }  \label{01}
\\
c_{1}\partial _{t}p_{1}+\alpha _{1}\mathrm{div}\left( \partial _{t}\mathbf{u}%
\right) -K_{1}\Delta p_{1}+g(p_{1}-p_{2}) &=&h_{1}\text{, }  \label{02} \\
c_{1}\partial _{t}p_{2}+\alpha _{2}\mathrm{div}\left( \partial _{t}\mathbf{u}%
\right) -K_{2}\Delta p_{2}-g(p_{1}-p_{2}) &=&h_{2}  \label{03}
\end{eqnarray}%
where $u$ is the displacement of the medium; $\lambda $ and $\mu $ are the
dilation and shear moduli of elasticity, respectively; $p_{i}$ is the
pressure of the fluid in phase $\left( i\right) $; $c_{i}$ the
compressibility, $K_{i}$ the permeability and $\alpha _{i}$ is the
Biot-Willis parameters \cite{bw2}. We note that if we let the volume of
fissures shrinks to zero so that $c_{2},\alpha _{2},K_{2},g$ become
negligible then the system \ (\ref{01})-(\ref{03}) reduces to the classical
Biot system with single porosity \cite{biot}:
\begin{eqnarray}
-\mu \Delta \mathbf{u}-(\lambda +\mu )\nabla (\mathrm{div}\mathbf{u})+\alpha
_{1}\nabla p_{1} &=&\mathbf{f},  \label{04} \\
c_{1}\partial _{t}p_{1}+\alpha _{1}\mathrm{div}\left( \partial _{t}\mathbf{u}%
\right) -K_{1}\Delta p_{1} &=&h_{1}.  \label{05}
\end{eqnarray}%
On the other hand, by neglecting the deformation effects $\lambda ,$ $\mu $
and $\alpha _{i}$ the system (\ref{01})-(\ref{03}) reduces to the BZK model
\cite{bzk}:%
\begin{eqnarray}
c_{1}\partial _{t}p_{1}-K_{1}\Delta p_{1}+g\left( p_{1}-p_{2}\right) &=&h_{1}%
\text{,}  \label{06} \\
c_{2}\partial _{t}p_{2}-K_{2}\Delta p_{2}-g\left( p_{1}-p_{2}\right) &=&h_{2}
\label{07}
\end{eqnarray}%
\

Aifantis' theory of consolidation with the concept of double porosity unify
then the proposed models (\ref{04})-(\ref{05}) of Biot for consolidation of
deformable porous media with single porosity and (\ref{06})-(\ref{07}) of
BZK model for fluid flow through undeformable porous media with double
porosity. Note also that a mathematical justification of the Aifantis model
has been established in \cite{ain2}. More precisely, it is considered
micro-models with periodically distributed poroelastic inclusions, embedded
in an extra poroelastic matrix, with imperfect contact on the interface. The
micro-model is based on Biot's system for consolidation processes with
interfacial barrier formulation. The macro-model is then derived by means of
the two-scale convergence technique and it reads as follows:%
\begin{gather}
-\mathrm{div}\sigma \left( \mathbf{u}\right) +\alpha _{1}\nabla p_{1}+\alpha
_{2}\nabla p_{2}=\mathbf{f},  \label{HP1} \\
\partial _{t}\left( \widetilde{c_{1}}p_{1}+\beta _{1}:e\left( \mathbf{u}%
\right) \right) -\mathrm{div}\left( K_{1}\nabla p_{1}\right) +\widetilde{g}%
\left( p_{1}-p_{2}\right) =h_{1}  \label{HP12} \\
\partial _{t}\left( \widetilde{c_{2}}p_{2}+\beta _{2}:e\left( \mathbf{u}%
\right) \right) -\mathrm{div}\left( K_{2}\nabla p_{2}\right) -\widetilde{g}%
\left( p_{1}-p_{2}\right) =h_{2}  \label{HP2}
\end{gather}%
where $\sigma ,$ $\alpha _{i},\beta _{i}$ and $K_{i}$ are some effective
tensors, $i=1,2$. See \cite{ain2} for more details. It is then worth
pointing out that the Aifantis model (\ref{01})-(\ref{03}) can be seen as a
special case of the homogenized model (\ref{HP1})-(\ref{HP2}) ($\beta
_{i}=\alpha _{i}=\gamma _{i}\mathrm{I}_{3},$ $\gamma _{i}$ being a scalar
and $\mathrm{I}_{3}$ the identity matrix).

In this paper, we consider a family of microscopic models for the fluid flow
in a periodic poroelastic medium made of two constituents : the matrix and
the inclusions, where the material properties change rapidly on a small
scale characterized by a parameter $\varepsilon $\ representing the
periodicity of the medium. We shall make the essential assumption that these
inclusions\ have sizes large enough compared with the sizes of pores so that
it makes sense to consider these media as poroelastic materials.

An interesting question is to investigate the limiting behavior of such
media when the flow in the inclusions presents very high frequency spatial
variations as a result of a relatively very low permeability when comparing
to the matrix permeability, since pore flow velocities in the porous matrix
can be high compared to movement through the interconnected pore spaces in
the inclusions. This leads especially to rescale the flow\ potential in the
inclusions by $\varepsilon ^{2}$, as in Arbogast \& \textit{al}. \cite{adh1}%
. The main objective of this paper is to derive a general model from the
point of view of homogenization theory. It will be seen that the macro-model
is in some sense the limit of a family of periodic micro-models in which the
size of the periodicity approach zero. It is shown that the overall behavior
of fluid flow in such heterogeneous media with low permeability at the
micro-scale may present memory terms. It is also shown that in such
anisotropic media, with different coupling interaction properties in
different directions, the Biot-Willis parameters are, as in \cite{ain2},
matrices and no longer scalars, as it is usually considered in the
poroelasticity litterature, since it is assumed there that the medium is
homogeneous and isotropic.

The paper is organized as follows. In the next section \ref{s2}, we give the
geometrical setting, the family of the periodic micro-models, and state the
main result of the paper. Section \ref{s3} is devoted to the proof of the
main result with the help of the two-scale convergence technique. Finally,
in section \ref{s4}, we conclude the paper with some remarks.

\section{Setting of the micro-model and main result\label{s2}}

The aim of this section is to provide a detailed set up of the studied
microstructure problem, introduce some necessary notations, basic
mathematical tools as well as the notion of two-scale convergence, auxiliary
problems, and then formulate the main result of the paper.

We consider $\Omega $\ a bounded and smooth domain in $\mathbb{R}^{3}$, $%
\varepsilon >0$ a sufficiently small parameter ($\varepsilon \ll 1$) and $%
Y=]0,1[^{3}$\ the generic cell of periodicity. We assume that $Y$\ is
divided as $Y=Y_{1}\cup Y_{2}\cup \Gamma $\ where $Y_{1},$\ $Y_{2}$\ are two
connected open subsets of $Y$ and $\Gamma $ is a smooth surface that
separates them. They are such that $\overline{Y_{2}}\subset Y$,$\ Y_{1}\cap
Y_{2}=\emptyset $,\ $\Gamma =\overline{Y_{1}}\cap \overline{Y_{2}}=\partial
Y_{2}$ and $\partial Y_{1}=\Gamma \cup \partial Y$. We denote $\ \mathbf{n}%
=\left( n_{i}\right) _{1\leq i\leq 3}$ the unit normal vector on $\Gamma $
pointing outward with respect to $Y_{1}$. Let $\chi _{1},\chi _{2}$ denote
respectively the characteristic function of $Y_{1}$, $Y_{2}$ extended by $Y$%
-periodicity to $\mathbb{R}^{3}$. Denote for $x\in \Omega $, $\chi
_{i}^{\varepsilon }(x)=\chi _{i}(x/\varepsilon )$ and set \textbf{\ }%
\begin{equation*}
\Omega _{i}^{\varepsilon }=\{x\in \Omega :\chi _{i}^{\varepsilon }(x)=1\}\
\text{and }\ \Gamma ^{\varepsilon }=\overline{\Omega _{1}^{\varepsilon }}%
\cap \overline{\Omega _{2}^{\varepsilon }}\text{.}
\end{equation*}%
Let $Z_{i}=\cup _{e\in \mathbb{Z}^{3}}\left( Y_{i}+e\right) $, $i=1,2$. As
in \cite{all}, we shall assume that\ the subset $Z_{1}$ is smooth and
connected open subset of $\mathbb{R}^{3}$.

With the above assumptions, the material occupying the domain $\Omega
_{2}^{\varepsilon }$\ is then embedded in the material occupying $\Omega
_{1}^{\varepsilon }$, and the interface $\Gamma ^{\varepsilon }$ is the
boundary $\partial \Omega _{2}^{\varepsilon }$. We observe that the boundary
of $\Omega _{1}^{\varepsilon }$ consists of two parts the outer boundary $%
\partial \Omega $\ and $\Gamma ^{\varepsilon }$.\textbf{\ }Usually, the
material $\Omega _{1}^{\varepsilon }$\ is referred to the matrix material
while the material $\Omega _{2}^{\varepsilon }$\ to the inclusions. Note
that no connectedness assumption is made on the material part $\Omega
_{2}^{\varepsilon }$.

Let $T>0$ and $t\in \left[ 0,T\right] $ denote the time variable. We set the
space-time domains $Q=\left( 0,T\right) \times \Omega $, $\Sigma =\left(
0,T\right) \times \Gamma $, $Q_{i}^{\varepsilon }=\left( 0,T\right) \times
\Omega _{i}^{\varepsilon }$, and $\Sigma ^{\varepsilon }=\left( 0,T\right)
\times \Gamma ^{\varepsilon }$.

Let us assume that each phase ($\Omega _{1}^{\varepsilon }$, $\Omega
_{2}^{\varepsilon }$) is occupied by a porous and deformable material
through which a slightly compressible and viscous fluid flow diffuses. Let $%
\mathbf{u}_{i}^{\varepsilon }$ denote the displacement of the medium $\Omega
_{i}^{\varepsilon }$, $i=1,2$. The equation of motion in $\Omega
_{1}^{\varepsilon }\cup \Omega _{2}^{\varepsilon }$ is given by
\begin{eqnarray}
-\mathrm{div}\sigma _{1}^{\varepsilon } &=&\mathbf{f}_{1},\ \text{in }\Omega
_{1}^{\varepsilon }\text{,}  \label{3} \\
-\mathrm{div}\sigma _{2}^{\varepsilon } &=&\mathbf{f}_{2},\ \text{in }\Omega
_{2}^{\varepsilon }  \label{3-0}
\end{eqnarray}%
where $\sigma _{i}^{\varepsilon }$ is the stress tensor which satisfies a
constitutive equation of linear poroelasticity of the form \cite{cou}:
\begin{equation}
\sigma _{i}^{\varepsilon }=\mathbb{A}_{i}^{\varepsilon }\mathrm{e}\left(
\mathbf{u}_{1}^{\varepsilon }\right) -\alpha _{i}^{\varepsilon
}p_{i}^{\varepsilon }\mathrm{I}_{3},\ \text{in }\Omega _{i}^{\varepsilon }
\label{3-1}
\end{equation}%
and $\mathbf{f}_{i\,}\in L^{2}\left( \Omega \right) ^{3}$ is the volume
distributed force in the corresponding medium, $i=1,2$. It is assumed that $%
\mathbf{f}_{i\,}$ is independent of $\varepsilon $. In (\ref{3-1}), $\mathbb{%
A}_{1}^{\varepsilon }$ and $\mathbb{A}_{2}^{\varepsilon }$ are fourth rank
elasticity tensors, $\mathrm{e}\left( \mathbf{\cdot }\right) $ is the
linearized strain tensor, $\mathrm{I}_{3}$ is the identity matrix, $%
p_{i}^{\varepsilon }$ is the pressure and $\alpha _{i}^{\varepsilon }$ is
the Biot-Willis parameter in the poroelastic material $\Omega
_{i}^{\varepsilon }$ \cite{bw2}.

Let $c_{1}^{\varepsilon }$\ (resp. $c_{2}^{\varepsilon }$) and $%
K_{1}^{\varepsilon }$ (resp. $K_{2}^{\varepsilon }$) denote respectively the
porosity and the permeability of the medium $\Omega _{1}^{\varepsilon }$
(resp. $\Omega _{2}^{\varepsilon }$). The equation for mass conservation in
each phase reads as follows:
\begin{eqnarray}
&&\partial _{t}\left( c_{1}^{\varepsilon }p_{1}^{\varepsilon }+\alpha
_{1}^{\varepsilon }\mathrm{div}\mathbf{u}_{1}^{\varepsilon }\right) -\mathrm{%
div}\left( K_{1}^{\varepsilon }\nabla p_{1}^{\varepsilon }\right) =0,\text{
in }\Omega _{1}^{\varepsilon },  \label{4} \\
&&\partial _{t}\left( c_{2}^{\varepsilon }p_{2}^{\varepsilon }+\alpha
_{2}^{\varepsilon }\mathrm{div}\mathbf{u}_{2}^{\varepsilon }\right) -\mathrm{%
div}\left( K_{2}^{\varepsilon }\nabla p_{2}^{\varepsilon }\right) =0\text{
in }\Omega _{2}^{\varepsilon }\text{.}  \label{5}
\end{eqnarray}%
On the interface $\Gamma ^{\varepsilon }$, we associate to (\ref{3})-(\ref%
{3-0}) the following transmission conditions:
\begin{equation}
\mathbf{u}_{1}^{\varepsilon }=\mathbf{u}_{2}^{\varepsilon }\text{, }\sigma
_{1}^{\varepsilon }\cdot \mathbf{n}^{\varepsilon }=\sigma _{2}^{\varepsilon
}\cdot \mathbf{n}^{\varepsilon }  \label{7}
\end{equation}%
and to (\ref{4})-(\ref{5}) the well-known open-pore conditions:%
\begin{equation}
\left( K_{1}^{\varepsilon }\nabla p_{1}^{\varepsilon }\right) \cdot \mathbf{n%
}^{\varepsilon }=\left( K_{2}^{\varepsilon }\nabla p_{2}^{\varepsilon
}\right) \cdot \mathbf{n}^{\varepsilon },\ \left( K_{1}^{\varepsilon }\nabla
p_{1}^{\varepsilon }\right) \cdot \mathbf{n}^{\varepsilon }=-g^{\varepsilon
}\left( p_{1}^{\varepsilon }-p_{2}^{\varepsilon }\right) .\text{ }  \label{9}
\end{equation}%
where $\mathbf{n}^{\varepsilon }$ stands for the unit normal vector on $%
\Gamma ^{\varepsilon }$ pointing outward with respect to $\Omega
_{1}^{\varepsilon }$, and $g^{\varepsilon }$ is the hydraulic permeability
of the thin layer $\Gamma ^{\varepsilon }$. The interface conditions (\ref{9}%
) are also known as \textsl{Deresiewicz-Skalak boundary conditions} \cite{ds}%
. Taking the limit on the thickness of the thin layer, one can prove that
these conditions are the only ones that are fully consistent with the
validity of the poroelasticity's equations throughout heterogeneous media
presenting interfaces across which the diffusion is discontinuous, see \cite%
{gur}. Observe that when $g^{\varepsilon }=\infty $, (\ref{9}) reduces to
the standard transmission condition, that is a perfect hydraulic contact on
the interface, and when $g^{\varepsilon }=0$, condition (\ref{9}) implies no
motion of the fluid relative to the solid. Here, in this paper we shall
assume that none of these conditions are fulfilled. See assumption (H4)
below.

On the exterior boundary $\partial \Omega \backslash \Gamma ^{\varepsilon }$%
, we assume the homogeneous Dirichlet boundary condition:%
\begin{equation}
\mathbf{u}_{1}^{\varepsilon }=\mathbf{0}\text{ and }p_{1}^{\varepsilon }=0%
\text{.}  \label{10}
\end{equation}

Finally, the system (\ref{4})-(\ref{10}) is supplemented by the following
initial conditions:%
\begin{eqnarray}
\mathbf{u}_{1}^{\varepsilon }\left( 0,\cdot \right) &=&\mathbf{0},\
p_{1}^{\varepsilon }\left( 0,\cdot \right) =0\text{ in }\Omega
_{1}^{\varepsilon }\text{,}  \label{11} \\
\mathbf{u}_{2}^{\varepsilon }\left( 0,\cdot \right) &=&\mathbf{0},\
p_{2}^{\varepsilon }\left( 0,\cdot \right) =0\text{ in }\Omega
_{2}^{\varepsilon }\text{.}  \label{12}
\end{eqnarray}

To deal with periodic homogenization with microstructures, we shall assume
the followings:

\begin{enumerate}
\item[(H1) ] There exists $Y$-periodic, fourth rank tensor-valued functions $%
\mathbb{A}_{i}\left( y\right) $, $i=1,2$ and continuous on $\mathbb{R}^{3}$
such that
\begin{equation*}
\mathbb{A}_{i}^{\varepsilon }\left( x\right) =\mathbb{A}_{i}\left( \frac{x}{%
\varepsilon }\right) ,\ \text{\ }x\in \Omega ,
\end{equation*}%
and
\begin{equation*}
\left( \mathbb{A}_{i}\left( y\right) \Xi :\Xi \right) \geq C\left( \Xi :\Xi
\right) \text{.}
\end{equation*}%
for all $y\in \mathbb{R}^{3}$\ and$\ \Xi \in \mathcal{M}_{\mathrm{sym}%
}^{3\times 3}\left( \mathbb{R}\right) $;

\item[(H2) ] There exist $Y$-periodic real-valued functions $c_{i}\left(
y\right) $, $i=1,2$ and continuous on $\mathbb{R}^{3}$ such that
\begin{equation*}
c_{i}^{\varepsilon }\left( x\right) =c_{i}\left( \frac{x}{\varepsilon }%
\right) ,\ x\in \Omega
\end{equation*}%
and $c_{i}\left( y\right) \geq C>0$ for all $y\in \mathbb{R}^{3}$;

\item[(H3) ] There exist $Y$-periodic matrix-valued functions $K_{i}\left(
y\right) $, $i=1,2$, continuous on $\mathbb{R}^{3}$ such that
\begin{equation}
K_{1}^{\varepsilon }\left( x\right) =K_{1}\left( \frac{x}{\varepsilon }%
\right) ,\ \ K_{2}^{\varepsilon }\left( x\right) =\varepsilon
^{2}K_{2}\left( \frac{x}{\varepsilon }\right) ,\text{\ \ \ }x\in \Omega
\label{sc1}
\end{equation}%
and
\begin{equation*}
\left\langle K_{i}\xi ,\xi \right\rangle \geq C\left\vert \xi \right\vert
^{2},\ i=1,2
\end{equation*}%
$\ $\ for all $y\in \mathbb{R}^{3}\ $and$\ \xi \in \mathbb{R}^{3}$;

\item[(H4) ] There exists a function $g\in \mathcal{C}\left( \mathbb{R}%
^{3}\right) ,\ Y$-periodic such that%
\begin{equation*}
g^{\varepsilon }\left( x\right) =\varepsilon g\left( x/\varepsilon \right)
,\ \ \ x\in \mathbb{R}^{3}\text{ and }\inf_{y\in \mathbb{R}^{3}}g\left(
y\right) \geq C>0.
\end{equation*}

\item[(H5) ] The Biot-Willis parameter $\alpha _{i}^{\varepsilon }$ is
defined a.e. in $\Omega $ as follows:
\begin{equation}
\alpha _{1}^{\varepsilon }\left( x\right) =\alpha _{1}\text{ for }x\in
\Omega _{1}^{\varepsilon }\text{ and }\alpha _{2}^{\varepsilon }\left(
x\right) =\varepsilon \alpha _{2}\text{ for }x\in \Omega _{2}^{\varepsilon }
\label{sc2}
\end{equation}%
where $\alpha _{i}$ is a positive constant, $i=1,2$.
\end{enumerate}

Here and thoughout this paper, the quantity $C$ denotes various positive
constants independent of $\varepsilon >0$, of the subscript $i=1,2$ and the
microscopic variable $y\in \mathbb{R}^{3}$.

\begin{remark}
We have chosen a particular scaling of the permeability coefficients in (\ref%
{sc1}). This means that the permeability is much larger in the network of
pores of the inclusions than in the porous matrix rocks. This gives that
both terms $\int_{\Omega _{1}^{\varepsilon }}\left\vert \nabla
p_{1}^{\varepsilon }\right\vert ^{2}dx$ and $\varepsilon ^{2}\int_{\Omega
_{2}^{\varepsilon }}\left\vert \nabla p_{2}^{\varepsilon }\right\vert ^{2}dx$
have the same order of magnitude and thus leading to a balance in potential
energies. For more details, we refer the reader to Arbogast, Douglas, and
Hornung \cite{adh1} (see also Allaire \cite{all}). In the same way, we also
have taken a special scaling factor of the Biot-Willis parameters in (\ref%
{sc2})leading to a balance in compressibility/dilation between the matrix
and inclusions.
\end{remark}

To set the mathematical framework of our Porblem, we need to introduce the
following spaces:
\begin{eqnarray*}
&&\mathbf{H}=H_{0}^{1}\left( \Omega \right) ^{3},\ \ \ \ \ \ \ \ \ \
L^{\varepsilon }=L^{2}\left( \Omega _{1}^{\varepsilon }\right) \times
L^{2}\left( \Omega _{2}^{\varepsilon }\right) , \\
&&\mathcal{E}_{1}^{\varepsilon }=\left\{ q\in H^{1}\left( \Omega
_{1}^{\varepsilon }\right) ;\ q_{|\Gamma }=0\right\} ,\ \mathcal{E}%
_{2}^{\varepsilon }=H^{1}\left( \Omega _{2}^{\varepsilon }\right) ,\
\mathcal{E}^{\varepsilon }=\mathcal{E}_{1}^{\varepsilon }\times \mathcal{E}%
_{2}^{\varepsilon }\text{.}
\end{eqnarray*}%
The space $\mathbf{H}$ is equipped with the standard norm: $||\mathbf{v}||_{%
\mathbf{H}}=||\mathrm{e}\left( \mathbf{v}\right) ||_{L^{2}\left( \Omega
\right) ^{3\times 3}}$ and $\mathcal{E}^{\varepsilon }$ with
\begin{equation*}
\left\Vert \left( q_{1},q_{2}\right) \right\Vert _{\mathcal{E}^{\varepsilon
}}^{2}=\left\Vert \nabla q_{1}\right\Vert _{L^{2}\left( \Omega
_{1}^{\varepsilon }\right) }^{2}+\varepsilon ^{2}\left\Vert \nabla
q_{2}\right\Vert _{L^{2}\left( \Omega _{2}^{\varepsilon }\right)
}^{2}+\varepsilon \left\Vert q_{1}-q_{2}\right\Vert _{L^{2}(\Gamma
^{\varepsilon })}^{2}\text{.}
\end{equation*}%
See S. Monsurro \cite{mons}. For a.e. $\left( t,x\right) \in Q$, we denote
\begin{eqnarray*}
\mathbf{u}^{\varepsilon }\left( t,x\right) &=&\chi _{1}^{\varepsilon }\left(
x\right) \mathbf{u}_{1}^{\varepsilon }\left( t,x\right) +\chi
_{2}^{\varepsilon }\left( x\right) \mathbf{u}_{2}^{\varepsilon }\left(
t,x\right) , \\
\ \ \ \mathbb{A}^{\varepsilon }\left( x\right) &=&\chi _{1}^{\varepsilon
}\left( x\right) \mathbb{A}_{1}^{\varepsilon }\left( x\right) +\chi
_{2}^{\varepsilon }\left( x\right) \mathbb{A}_{2}^{\varepsilon }\left(
x\right) , \\
\mathbf{f}^{\varepsilon }\left( x\right) &=&\chi _{1}^{\varepsilon }\left(
x\right) \mathbf{f}_{1}\left( x\right) +\chi _{2}^{\varepsilon }\left(
x\right) \mathbf{f}_{2}\left( x\right) \text{.}
\end{eqnarray*}%
Note that, thanks to the transmission condition (\ref{7}), the displacement $%
\mathbf{u}^{\varepsilon }\left( t,\cdot \right) $ lies in $\mathbf{H}$ for
a.e. $t\in \left( 0,T\right) $.

Throughout the paper the following notation will be used: if $\mathcal{F}$
is any functional space then $L_{T}^{p}\left( \mathcal{F}\right) $ denotes
the Bochner vector-valued functions space defined by $L_{T}^{p}\left(
\mathcal{F}\right) =L^{p}\left( 0,T;\mathcal{F}\right) $

The weak formulation of (\ref{4})-(\ref{12}) can be read as follows\textbf{:}

Find $\ \left( \mathbf{u}^{\varepsilon },p^{\varepsilon }\right) \in
L_{T}^{\infty }\left( \mathbf{H}\right) \times L_{T}^{2}\left( \mathcal{E}%
^{\varepsilon }\right) $, such that $p^{\varepsilon }=\left(
p_{1}^{\varepsilon },p_{2}^{\varepsilon }\right) \in L_{T}^{\infty }\left(
L^{\varepsilon }\right), $ 
\begin{equation*}
\partial _{t}\left( c_{1}^{\varepsilon }p_{1}^{\varepsilon }+\alpha _{1}%
\mathrm{div}\mathbf{u}^{\varepsilon }\right) \in L_{T}^{2}\left( \mathcal{E}%
_{1}^{\varepsilon }{}^{\ast }\right) ,\partial _{t}\left( c_{2}^{\varepsilon
}p_{2}^{\varepsilon }+\varepsilon \alpha _{2}\mathrm{div}\mathbf{u}%
^{\varepsilon }\right) \in L_{T}^{2}\left( \mathcal{E}_{2}^{\varepsilon
}{}^{\ast }\right)
\end{equation*}
 and for all $\mathbf{v}\in
\mathbf{H}$, $\left( q_{1},q_{2}\right) \in \mathcal{E}^{\varepsilon }$, we
have%

\begin{equation}
\int_{\Omega }\mathbb{A}^{\varepsilon }\mathrm{e}\left( \mathbf{u}%
^{\varepsilon }\right) \mathrm{e}\left( \mathbf{v}\right) dx+\int_{\Omega
_{1}^{\varepsilon }}\alpha _{1}\nabla p_{1}^{\varepsilon }\mathbf{v}%
dx+\int_{\Omega _{2}^{\varepsilon }}\alpha _{2}^{\varepsilon }\nabla
p_{2}^{\varepsilon }\mathbf{v}dx=\int_{\Omega }\mathbf{f}^{\varepsilon }%
\mathbf{v}dx\text{,}  \label{23}
\end{equation}%
\begin{eqnarray}
&&\left\langle \partial _{t}\left( c_{1}^{\varepsilon }p_{1}^{\varepsilon
}+\alpha _{1}\mathrm{div}\mathbf{u}^{\varepsilon }\right)
,q_{1}\right\rangle _{\mathcal{E}_{1}^{\varepsilon }{}^{\ast },\mathcal{E}%
_{1}^{\varepsilon }}+\int_{\Omega _{1}^{\varepsilon }}K_{1}^{\varepsilon
}\nabla p_{1}^{\varepsilon }\nabla q_{1}dx+  \notag \\
&&\left\langle \partial _{t}\left( c_{2}^{\varepsilon }p_{2}^{\varepsilon
}+\varepsilon \alpha _{2}\mathrm{div}\mathbf{u}^{\varepsilon }\right)
,q_{2}\right\rangle _{\mathcal{E}_{2}^{\varepsilon }{}^{\ast },\mathcal{E}%
_{2}^{\varepsilon }}+\int_{\Omega _{2}^{\varepsilon }}K_{2}^{\varepsilon
}\nabla p_{2}^{\varepsilon }\nabla q_{2}dx+  \notag \\
&&\int_{\Gamma ^{\varepsilon }}g^{\varepsilon }\left( p_{1}^{\varepsilon
}-p_{2}^{\varepsilon }\right) \left( q_{1}-q_{2}\right) ds^{\varepsilon
}\left( x\right) =0,  \label{24} \\
&&\mathbf{u}^{\varepsilon }\left( 0,\cdot \right) =\mathbf{0}\text{,\ }\chi
_{1}\left( \cdot \right) p_{1}^{\varepsilon }\left( 0,\cdot \right) +\chi
_{2}\left( \cdot \right) p_{2}^{\varepsilon }\left( 0,\cdot \right) =0\text{
a.e. in }\Omega \text{.}  \label{241}
\end{eqnarray}%
\textbf{\ \ } Here and throughout this paper $dx$ and $ds^{\varepsilon }\left( x\right) $ stand
respectively for the Lebesgue measure on $\mathbb{R}^{3}$ and the Hausdorff
measure on $\Gamma ^{\varepsilon }$.

Using assumptions (H1)-(H5), we can establish the following existence and
uniqueness result whose proof is a slight modification of that\ given by R.
E. Showalter and B. Momken \cite{showmomk} and therefore will be omitted.

\begin{theorem}
Assume\ that (H1)-(H5) hold. Then, for any sufficiently small $\varepsilon
>0 $ and $\mathbf{f}^{\varepsilon }\in \mathbf{L}^{2}\left( \Omega \right) $%
, there exists a unique couple $\left( \mathbf{u}^{\varepsilon
},p^{\varepsilon }\right) \in L_{T}^{\infty }\left( \mathbf{H}\right) \times
L_{T}^{2}\left( \mathcal{E}^{\varepsilon }\right) $, solution of the weak
system (\ref{23})-(\ref{241}), such that
\begin{equation}
\Vert \mathbf{u}^{\varepsilon }\Vert _{L_{T}^{\infty }\left( \mathbf{H}%
\right) }+\Vert p^{\varepsilon }\Vert _{L_{T}^{2}\left( \mathcal{E}%
^{\varepsilon }\right) }+\Vert p^{\varepsilon }\Vert _{L_{T}^{\infty }\left(
L^{\varepsilon }\right) }\leq C\text{.}  \label{30}
\end{equation}
\end{theorem}

Now, thanks to the a priori estimates (\ref{30}), one is led to study the
limiting behavior of\ the sequence $\left( \mathbf{u}^{\varepsilon
},p^{\varepsilon }\right) $ as $\varepsilon $ approaches $0$. To do this, we
shall use the two-scale convergence technique that we shall recall hereafter.

First, we define $C_{\#}(Y)$ to be the space of all continuous functions on $%
\mathbb{R}^{3}$ which are $Y$-periodic. Let the space $L_{\#}^{2}\left(
Y\right) $ (resp. $L_{\#}^{2}\left( Y_{i}\right) $, $i=1,2$) to be all
functions belonging to $L_{\mathrm{loc}}^{2}\left( \mathbb{R}^{3}\right) $
(resp. $L_{\mathrm{loc}}^{2}\left( Z_{i}\right) $) which are $Y$-periodic,
and $H_{\#}^{1}\left( Y\right) $ (resp. $H_{\#}^{1}\left( Y_{i}\right) $) to
be the space of those functions together with their derivatives belonging to
$L_{\#}^{2}\left( Y\right) $ (resp. $L_{\#}^{2}\left( Z_{i}\right) $).

Now, we recall the definition and main results concerning the method of
two-scale convergence. For more details, we refer the reader to \cite{all,
adh, ngue}.

\begin{definition}
\label{def1}A sequence $\left( v^{\varepsilon }\right) \ $in $L^{2}\left(
\Omega \right) $ two-scale converges to $v\in L^{2}\left( \Omega \times
Y\right) $ (we write $v^{\varepsilon }\overset{2-s}{\rightharpoonup }v$) if,
for any admissible test function $\varphi \in L^{2}\left( \Omega ;\mathcal{C}%
_{\#}(Y)\right) $,
\begin{equation*}
\lim_{\varepsilon \rightarrow 0}\int_{\Omega }v^{\varepsilon }\left(
x\right) \varphi \left( x,\frac{x}{\varepsilon }\right) dx=\int_{\Omega
\times Y}v\left( x,y\right) \varphi \left( x,y\right) dxdy\text{.}
\end{equation*}
\end{definition}

\begin{theorem}
\label{t1}Let $(v^{\varepsilon })$ be a sequence of functions in $%
L^{2}(\Omega )$ which is uniformly bounded. Then, there exist $v\in
L^{2}(\Omega \times Y)$ and a subsequence of $(v^{\varepsilon })$ which
two-scale converges\ to $v$.
\end{theorem}

\begin{theorem}
\label{t2}Let $(v^{\varepsilon })$ be a uniformly bounded sequence in $%
H^{1}(\Omega )$ (resp. $H_{0}^{1}(\Omega )$). Then there exist $v\in
H^{1}\left( \Omega \right) $ (resp. $H_{0}^{1}(\Omega )$) and $\hat{v}\in
L^{2}(\Omega ;H_{\#}^{1}(Y)/\mathbb{R})$ such that, up to a subsequence,%
\begin{equation*}
v^{\varepsilon }\overset{2-s}{\rightharpoonup }v;\ \ \ \ \ \nabla
v^{\varepsilon }\overset{2-s}{\rightharpoonup }\nabla v+\nabla _{y}\hat{v}.
\end{equation*}
\end{theorem}

Here and in the sequel the subscript $y$ on a differential operator denotes
the derivative with respect to $y$.

\begin{theorem}
\label{t3}Let $(v^{\varepsilon })$ be a sequence of functions in $%
H^{1}(\Omega )$ such that
\begin{equation*}
\left\Vert v^{\varepsilon }\right\Vert _{L^{2}\left( \Omega \right)
}+\varepsilon \left\Vert \nabla v^{\varepsilon }\right\Vert _{L^{2}\left(
\Omega \right) ^{3}}\leq C\text{.}
\end{equation*}%
Then, there exist $v\in L^{2}\left( \Omega ;H_{\#}^{1}(Y)\right) $ and a
subsequence of $\left( v^{\varepsilon }\right) $, still denoted by $\left(
v^{\varepsilon }\right) $ such that
\begin{equation*}
v^{\varepsilon }\overset{2-s}{\rightharpoonup }v,\ \ \ \ \ \varepsilon
\nabla v^{\varepsilon }\overset{2-s}{\rightharpoonup }\nabla _{y}v
\end{equation*}%
and for every $\varphi \in \mathcal{D}\left( \Omega ;\mathcal{C}%
_{\#}(Y)\right) $, we have:%
\begin{equation*}
\lim_{\varepsilon \rightarrow 0}\int_{\Gamma ^{\varepsilon }}v^{\varepsilon
}\left( x\right) \varphi \left( x,\frac{x}{\varepsilon }\right)
ds^{\varepsilon }\left( x\right) =\int_{\Omega \times \Gamma }v\left(
x,y\right) \varphi \left( x,y\right) dxds\left( y\right) .
\end{equation*}%
Here and in the sequel $ds\left( y\right) $ denotes the Hausdorff measure on
$\Gamma $.
\end{theorem}

The notion of two-scale convergence can easily be generalized to
time-dependent functions without affecting the results stated above.
According to \cite{cs}, we give the following:

\begin{definition}
\label{def2}We say that a sequence $\left( v^{\varepsilon }\right) $ in $%
L^{2}\left( Q\right) $ two-scale converges to $v\in L^{2}\left( Q\times
Y\right) $ (we always write $v^{\varepsilon }\overset{2-s}{\rightharpoonup }%
v $) if, for any $\varphi \in L^{2}\left( Q;\mathcal{C}_{\#}(Y)\right) $:%
\begin{equation*}
\lim_{\varepsilon \rightarrow 0}\int_{Q}v^{\varepsilon }\left( t,x\right)
\varphi \left( t,x,\frac{x}{\varepsilon }\right) dtdx=\int_{Q\times
Y}v\left( t,x,y\right) \varphi \left( t,x,y\right) dtdxdy\text{.}
\end{equation*}
\end{definition}

\begin{remark}
\label{rem1}The results stated above still hold for the case of
time-dependent sequences. For if $\left( v^{\varepsilon }\right) $ is a
uniformly bounded sequence in $L^{2}\left( Q\right) $, there exists then $%
v\in L^{2}\left( Q\right) $ such that, up to a subsequence, $v^{\varepsilon }%
\overset{2-s}{\rightharpoonup }v$ in the sense of Def. \ref{def2}. Moreover,
if $\left( v^{\varepsilon }\right) $ is uniformly bounded in $%
L_{T}^{2}\left( H^{1}\left( \Omega \right) \right) $, then up to a
subsequence, there exist $v\in L_{T}^{2}\left( H^{1}\left( \Omega \right)
\right) $ and $v_{0}\in L^{2}\left( Q;H_{\#}^{1}\left( Y\right) /\mathbb{R}%
\right) $ such that $v^{\varepsilon }\overset{2-s}{\rightharpoonup }v$ and $%
\nabla v^{\varepsilon }\overset{2-s}{\rightharpoonup }\nabla v+\nabla
_{y}v_{0}$. On the other hand, if a sequence $\left( v^{\varepsilon }\right)
$ is such that
\begin{equation*}
\left\Vert v^{\varepsilon }\right\Vert _{L^{2}\left( Q\right) }+\varepsilon
\left\Vert \nabla v^{\varepsilon }\right\Vert _{L^{2}\left( Q\right) }\leq C,
\end{equation*}%
then, up to a subsequence, there exists $v\in L_{T}^{2}\left(
H_{\#}^{1}\left( Y\right) \right) $ such that $v^{\varepsilon }\overset{2-s}{%
\rightharpoonup }v$ and $\varepsilon \nabla _{y}v^{\varepsilon }\overset{2-s}%
{\rightharpoonup }\nabla _{y}v$.
\end{remark}

In order to state the main result, we shall give in the sequel some
notations. Let us first introduce the three following auxilliary problems.
For $j,k\in \left\{ 1,2,3\right\} $, let$\ \mathbf{w}^{jk}\in \left(
H_{\#}^{1}\left( Y\right) /\mathbb{R}\right) ^{3}$ be the solution to the
following microscopic system:
\begin{equation*}
\left\{
\begin{array}{l}
-\mathrm{div}_{y}\left( \mathbb{A}_{1}\mathrm{e}_{y}\left( \mathbf{w}^{jk}+%
\mathbf{d}^{jk}\right) \right) =0\text{ a.e. in }Y_{1}\text{,} \\
-\mathrm{div}_{y}\left( \mathbb{A}_{2}\mathrm{e}_{y}\left( \mathbf{w}^{jk}+%
\mathbf{d}^{jk}\right) \right) =0\text{ a.e. in }Y_{2}\text{,} \\
\mathbb{A}_{1}\mathrm{e}_{y}\left( \mathbf{w}^{jk}+\mathbf{d}^{jk}\right)
\cdot \mathbf{n}=\mathbb{A}_{2}\mathrm{e}_{y}\left( \mathbf{w}^{jk}+\mathbf{d%
}^{jk}\right) \cdot \mathbf{n}\text{ a.e. on }\Gamma \text{,} \\
y\longmapsto \ \mathbf{w}^{jk}\ Y\text{-periodic}%
\end{array}%
\right.
\end{equation*}%
where $\mathbf{d}^{jk}\left( y\right) =\left( y_{j}\delta _{lk}\right)
_{1\leq l\leq 3}$ and $\left( \delta _{kj}\right) $ is the Kr\"{o}necker
symbol. For $j=1,2,3$, let $\pi _{j}\in H^{1}\left( Y_{1}\right) /\mathbb{R}$
be the solution of the following stationary micro-pressure equation:
\begin{equation*}
\left\{
\begin{array}{l}
-\mathrm{div}_{y}\left( K_{1}\left( \nabla \pi _{j}+e_{j}\right) \right) =0%
\text{ in }Y_{1}\text{,} \\
K_{1}\left( \nabla \pi _{j}+e_{j}\right) \cdot \mathbf{n}=0\text{ on }\Gamma
\text{,} \\
y\longmapsto \pi _{j}\text{ }Y\text{-periodic}%
\end{array}%
\right.
\end{equation*}%
where $e_{j}$ is the $j^{\text{th}}$ vector of the canonical basis of $%
\mathbb{R}^{3}$. Let $\zeta \in L_{T}^{\infty }\left( H_{\#}^{1}\left(
Y_{2}\right) \right) $ be the unique solution to the following non
micro-pressure problem of the Robin type:
\begin{equation*}
\left\{
\begin{array}{l}
\partial _{t}\left( c_{2}\zeta \right) -\mathrm{div}_{y}\left( K_{2}\nabla
_{y}\zeta \right) \text{%
$=$%
\ }0\text{ a.e.in }\left( 0,T\right) \times Y_{2}\text{,} \\
K_{2}\nabla _{y}\zeta \cdot \mathbf{n}\text{%
$=$%
\ }-g\left( y\right) \left[ 1-\zeta \right] \text{ a.e. on }\Sigma \text{,}
\\
y\longmapsto \zeta \text{ }Y\text{-periodic,} \\
\zeta \left( 0,y\right) =0\text{,\ a.e. }y\in Y_{2}\text{.}%
\end{array}%
\right.
\end{equation*}%
Now, let us define the homogenized fourth rank tensor $\widetilde{\mathbb{A}}%
=\left( \tilde{a}_{j_{1}j_{2}j_{3}j_{4}}\right) _{1\leq
j_{1},j_{2},j_{3},j_{4}\leq 3}$, where the coefficients are given by
\begin{equation*}
\tilde{a}_{j_{1}j_{2}j_{3}j_{4}}=\sum_{k_{1},k_{2}=1}^{3}%
\int_{Y}a_{j_{1}j_{2}k_{1}k_{2}}\left( y\right) \left( \delta
_{j_{1}k_{1}}\delta _{j_{2}k_{2}}+\mathrm{e}_{k_{1}k_{2},y}\left( \mathbf{w}%
^{j_{3}j_{4}}\right) \left( y\right) \right) dy\text{.}
\end{equation*}%
Here $\left( a_{jklm}\right) $ are the coefficients of the elasticity tensor
$\mathbb{A}$ which are given by
\begin{equation}
\mathbb{A}\left( y\right) =\chi _{1}\left( y\right) \mathbb{A}_{1}\left(
y\right) +\chi _{2}\left( y\right) \mathbb{A}_{2}\left( y\right)  \label{n0}
\end{equation}%
for a.e. $y\in Y$, and $\mathrm{e}_{jk,y}\left( \mathbf{\cdot }\right) $ is
the linearized elasticity strain tensor where the derivatives are taken with
respect to the microscopic variable $y$. Let also define the following
homogenized tensors:
\begin{equation}
\tilde{\sigma}(\mathbf{u})=\left( \tilde{\sigma}_{jk}(\mathbf{u})\right) ,\
\tilde{K}=\left( \tilde{K}_{jk}\right) ,\ B=\left( b_{jk}\right) ,\ \Lambda
=\left( \lambda _{jk}\right)  \label{n1}
\end{equation}%
where for $j,k\in \left\{ 1,2,3\right\} $
\begin{eqnarray}
&&\tilde{\sigma}_{jk}(\mathbf{u})=\sum_{l,m=1}^{3}\tilde{a}_{jklm}\mathrm{e}%
_{lm}(\mathbf{u})\text{,\ }  \label{n2} \\
&&\tilde{K}_{jk}=\int_{Y_{1}}K_{1}\left( y\right) \left( \nabla _{y}\pi
_{j}+e_{j}\right) \left( \nabla \pi _{k}+e_{k}\right) dy,  \label{n3} \\
&&b_{jk}=\alpha _{1}\left( \left\vert Y_{1}\right\vert \delta
_{jk}+\int_{\Gamma }\pi _{k}\left( y\right) n_{j}ds\left( y\right) \right) ,
\label{n4} \\
&&\lambda _{jk}=\alpha _{1}\int_{Y_{1}}\sum_{l=1}^{3}\left( \delta
_{jl}\delta _{kl}+\frac{\partial w_{l}^{jk}}{\partial y_{l}}\right) dy.
\label{n5}
\end{eqnarray}%
Here $\left\vert Y_{i}\right\vert $ denotes the volume of $Y_{i}$ and $%
\left( w_{l}^{ij}\right) _{1\leq l\leq 3}$ are the components of $\mathbf{w}%
^{ij}$. Finally let us define the following averaging quantities%
\begin{eqnarray}
\mathbf{f} &=&|Y_{1}|\mathbf{f}_{1}+|Y_{2}|\mathbf{f}_{2},  \label{n6} \\
\tilde{c} &=&\int_{Y_{1}}c_{1}\left( y\right) dy,  \label{n7} \\
\tilde{g} &=&\int_{\Gamma }g\left( y\right) ds\left( y\right)  \label{n8}
\end{eqnarray}%
and the time-dependent functions%
\begin{eqnarray}
\theta \left( t,\tau \right) &=&\alpha _{2}\int_{\Gamma }\partial _{t}\zeta
\left( t-\tau ,y\right) \mathbf{n}ds\left( y\right) ,\   \label{n10} \\
\eta \left( t,\tau \right) &=&-\int_{\Gamma }g\left( y\right) \partial
_{t}\zeta \left( t-\tau ,y\right) ds\left( y\right) \text{.}  \label{n11}
\end{eqnarray}

With these notations, we are now ready to give the main result of the paper:

\begin{theorem}
\label{thp}Let $\left( \mathbf{u}^{\varepsilon },p^{\varepsilon }\right) \in
L^{\infty }\left( 0,T;\mathbf{H}\right) \times L^{2}\left( 0,T;\mathcal{E}%
^{\varepsilon }\right) $ be the solution of the weak system (\ref{23}).
Then, up to a subsequence, there exists $\left( \mathbf{u},p\right) \in
L^{2}\left( 0,T;\mathbf{H}_{0}^{1}\left( \Omega \right) \times
H_{0}^{1}\left( \Omega \right) \right) $ such that
\begin{eqnarray*}
\mathbf{u}^{\varepsilon } &\rightharpoonup &\mathbf{u}\text{ in }L^{2}\left(
0,T;H_{0}^{1}\left( \Omega \right) \right) \text{ weakly,} \\
p_{1}^{\varepsilon } &\rightharpoonup &p_{1}\text{ in }L^{2}\left( Q\right)
\text{ weakly,} \\
p_{2}^{\varepsilon } &\rightharpoonup &\int_{Y_{2}}p_{2}\left( y\right) dy%
\text{ in }L^{2}\left( Q\right) \text{ weakly,}
\end{eqnarray*}%
where $p=\left( p_{1},\int_{Y_{2}}p_{2}\left( y\right) dy\right) ,$
\begin{equation*}
p_{2}\left( t,x,y\right) =\int_{0}^{t}p_{1}\left( \tau ,x\right) \partial
_{t}\zeta \left( t-\tau ,y\right) d\tau ,\ \text{a.e. }\left( t,x,y\right)
\in Q\times Y_{2}\text{.}
\end{equation*}%
and the couple $\left( \mathbf{u},p_{1}\right) $ is a solution to the
homogenized model:%
\begin{eqnarray*}
&&-\mathrm{div}\tilde{\sigma}(\mathbf{u})+B\nabla p_{1}+\int_{0}^{t}\theta
\left( t,\tau \right) p_{1}\left( \tau ,x\right) d\tau \text{%
$=$%
\ }\mathbf{f}\text{, a.e. in }Q, \\
&&\partial _{t}\left( \tilde{c}p_{1}+\Lambda :\mathrm{e}\left( \mathbf{u}%
\right) \right) -\mathrm{div}\left( \tilde{K}\nabla p_{1}\right) +\tilde{g}%
p_{1}-\int_{0}^{t}\eta \left( t,\tau \right) p_{1}\left( \tau ,x\right)
d\tau \text{%
$=$%
}0,\text{ a.e. in }Q, \\
&&\mathbf{u}\text{%
$=$%
\ }0\text{ and }\tilde{K}\nabla p_{1}\cdot \nu \text{%
$=$%
}0\text{ a.e. on }\Sigma , \\
&&\mathbf{u}\left( 0,x\right) \text{%
$=$%
\ }\mathbf{0}\text{\ a.e. in }\Omega \text{, }p_{1}\left( 0,x\right) \text{%
$=$%
\ }0\text{ a.e. in }\Omega \text{,}
\end{eqnarray*}

Here $\tilde{\sigma}$, $B$, $\theta $, $\mathbf{f}$\textbf{, } $\tilde{c}$, $%
\Lambda $, $\tilde{K}$, $\tilde{g}$ and $\eta $ are given in (\ref{n1})-(\ref%
{n11}).
\end{theorem}

\section{Proof of the main result\label{s3}}

As a direct application of the theorems listed above (Thms \ref{t1}-\ref{t3}%
) and the a priori estimates (\ref{30}), we give without proof the following
two-scale convergence result concerning the solutions $\left( \mathbf{u}%
^{\varepsilon },p^{\varepsilon }\right) $ of the Problem (\ref{23})-(\ref%
{241}).

\begin{theorem}
\label{t4}There exists a subsequence of $\left( \mathbf{u}^{\varepsilon
},p^{\varepsilon }\right) $, solution of (\ref{23})-(\ref{241}), still
denoted $\left( \mathbf{u}^{\varepsilon },p^{\varepsilon }\right) $, and
there exist
\begin{eqnarray*}
&&\mathbf{u}\in L_{T}^{\infty }\left( \mathbf{H}\right) \text{,\ \ \ }%
\mathbf{\hat{u}}\in L_{T}^{\infty }\left( L^{2}\left( \Omega
;H_{\#}^{1}\left( Y\right) /\mathbb{R}\right) \right) ^{3} \\
&&\  \\
&&p_{1}\in L_{T}^{\infty }\left( H_{0}^{1}\left( \Omega \right) \right)
\text{,\ \ }\ \hat{p}_{1}\in L^{2}(Q;H_{\#}^{1}\left( Y\right) /\mathbb{R})
\end{eqnarray*}%
and%
\begin{equation*}
p_{2}\in L_{T}^{\infty }\left( L^{2}\left( \Omega ;H_{\#}^{1}\left( Y\right)
\right) \right)
\end{equation*}%
such that, for a.e. $t\in \left( 0,T\right) $,
\begin{eqnarray}
&&\mathbf{u}^{\varepsilon }\left( t,\cdot \right) \overset{2-s}{%
\rightharpoonup }\mathbf{u}\left( t,\cdot \right) \mathbf{,}  \label{31} \\
&&\chi _{1}^{\varepsilon }p_{1}^{\varepsilon }\left( t,\cdot \right) \overset%
{2-s}{\rightharpoonup }\chi _{1}p_{1\,}\left( t,\cdot \right) ,  \label{32}
\\
&&\chi _{1}^{\varepsilon }p_{2}^{\varepsilon }\left( t,\cdot \right) \overset%
{2-s}{\rightharpoonup }\chi _{2}p_{2}\left( t,\cdot \right)  \label{321}
\end{eqnarray}%
in the sense of Def. \ref{def1}\ and
\begin{eqnarray}
&&\frac{\partial \mathbf{u}^{\varepsilon }}{\partial \mathbf{x}_{j}}\overset{%
2-s}{\rightharpoonup }\frac{\partial \mathbf{u}}{\partial \mathbf{x}_{j}}+%
\frac{\partial \mathbf{\hat{u}}}{\partial \mathbf{y}_{j}}\mathbf{,\ \ \ }%
j=1,2,3,  \label{33} \\
&&\chi _{1}^{\varepsilon }\nabla p_{1}^{\varepsilon }\overset{2-s}{%
\rightharpoonup }\chi _{1}\left( \nabla p_{1}+\nabla _{y}\hat{p}_{1}\right) ,
\label{34} \\
&&\varepsilon \chi _{2}^{\varepsilon }\nabla p_{2}^{\varepsilon }\overset{2-s%
}{\rightharpoonup }\chi _{2}\nabla _{y}p_{2}  \label{35}
\end{eqnarray}%
in the sense of Def. \ref{def2}. Moreover, the following convergence holds:
\
\begin{equation}
\lim_{\varepsilon \rightarrow 0}\int_{\Sigma ^{\varepsilon }}\varepsilon
\left( p_{1}^{\varepsilon }-p_{2}^{\varepsilon }\right) \psi ^{\varepsilon
}dtds^{\varepsilon }=\int_{Q\times \Gamma }\left( p_{1}-p_{2}\right) \psi
dtdxds,  \label{36}
\end{equation}%
for any $\psi \in \mathcal{D}\left( Q;\mathcal{C}_{\#}\left( Y\right)
\right) $ with $\psi ^{\varepsilon }\left( t,x\right) =\psi \left(
t,x,x/\varepsilon \right) $.
\end{theorem}

To determine the limiting equations of the system (\ref{23})-(\ref{241}), we
begin by choosing the adequate admissible test functions. Let $\mathbf{v}%
^{\varepsilon }\left( x\right) =\mathbf{v}\left( x\right) +\varepsilon
\mathbf{\hat{v}}\left( x,\dfrac{x}{\varepsilon }\right) $ where $\mathbf{v}%
\in \mathcal{D}\left( \Omega \right) ^{3}$ and $\mathbf{\hat{v}}\in \mathcal{%
D}\left( \Omega ;\mathcal{C}_{\#}^{\infty }\left( Y\right) \right) ^{3}$.
Let also $q_{1}^{\varepsilon }(t,x)=\varphi _{1}\left( t,x\right)
+\varepsilon \hat{\varphi}_{1}\left( t,x,\dfrac{x}{\varepsilon }\right) $
and $q_{2}^{\varepsilon }(t,x)=\varphi _{2}\left( t,x,\dfrac{x}{\varepsilon }%
\right) $ where $\varphi _{1}\in \mathcal{D}\left( \left( 0,T\right) \times
\bar{\Omega}\right) $ and $\varphi _{2}$, $\hat{\varphi}_{1}\in \mathcal{D}%
\left( Q;\mathcal{C}_{\#}^{\infty }\left( Y\right) \right) $. Taking $%
\mathbf{v}=\mathbf{v}^{\varepsilon }$ in (\ref{23}), we have:
\begin{eqnarray}
\int_{\Omega }\mathbf{f}^{\varepsilon }\mathbf{v}^{\varepsilon }dx
&=&\int_{\Omega }\mathbb{A}^{\varepsilon }\left( x\right) \mathrm{e}\left(
\mathbf{u}^{\varepsilon }\right) \mathrm{e}\left( \mathbf{v}^{\varepsilon
}\right) dx+\int_{\Omega _{1}^{\varepsilon }}\alpha _{1}\nabla
p_{1}^{\varepsilon }\mathbf{v}^{\varepsilon }dx+\varepsilon \int_{\Omega
_{2}^{\varepsilon }}\alpha _{2}\nabla p_{2}^{\varepsilon }\mathbf{v}%
^{\varepsilon }dx  \notag \\
&=&\int_{\Omega }\mathbb{A}^{\varepsilon }\left( x\right) \mathrm{e}\left(
\mathbf{u}^{\varepsilon }\right) \left( \mathrm{e}\left( \mathbf{v}\right)
\left( x\right) +\mathrm{e}_{y}\left( \mathbf{\hat{v}}\right) \left( x,\frac{%
x}{\varepsilon }\right) \right) dx+  \notag \\
&&\int_{\Omega }\left( \alpha _{1}\chi _{1}^{\varepsilon }\left( x\right)
\nabla p_{1}^{\varepsilon }+\varepsilon \alpha _{2}\chi _{2}^{\varepsilon
}\left( x\right) \nabla p_{2}^{\varepsilon }\right) \mathbf{v}\left(
x\right) dx+\varepsilon R_{1}^{\varepsilon },  \label{40}
\end{eqnarray}%
where
\begin{eqnarray*}
R_{1}^{\varepsilon } &=&\int_{\Omega }\mathbb{A}^{\varepsilon }\left(
x\right) \mathrm{e}\left( \mathbf{u}^{\varepsilon }\right) \mathrm{e}%
_{x}\left( \mathbf{w}\right) \left( x,\frac{x}{\varepsilon }\right)
dx+\alpha _{1}\int_{\Omega }\chi _{1}^{\varepsilon }\left( x\right) \nabla
p_{1}^{\varepsilon }\mathbf{w}\left( x,\frac{x}{\varepsilon }\right) dx \\
&&+\varepsilon \alpha _{2}\int_{\Omega }\chi _{2}^{\varepsilon }\left(
x\right) \nabla p_{2}^{\varepsilon }\mathbf{w}\left( x,\frac{x}{\varepsilon }%
\right) dx.
\end{eqnarray*}%
Observe that $R_{1}^{\varepsilon }=O\left( 1\right) $.

Now, we pass to the limit in (\ref{40}). In view of (\ref{33}), and since $%
\mathbb{A}^{t}\left( \mathrm{e}\left( \mathbf{v}\right) +\mathrm{e}%
_{y}\left( \mathbf{\hat{v}}\right) \right) $ is an admissible test function, the first integral in the l.h.s.
of (\ref{40}) converges to
\begin{equation}
\int_{\Omega \times Y}\mathbb{A}\left( \mathrm{e}\left( \mathbf{u}\right) +%
\mathrm{e}_{y}\left( \mathbf{\hat{u}}\right) \right) \left( \mathrm{e}\left(
\mathbf{v}\right) +\mathrm{e}_{y}\left( \mathbf{\hat{v}}\right) \right) dxdy
\label{401}
\end{equation}%
where the tensor $\mathbb{A}\left( y\right) $ is given by (\ref{n0}). In
view of Divergence Lemma and (\ref{34})-(\ref{35}), the second integral of
the l.h.s. of (\ref{40}) tends to
\begin{eqnarray}
&&\alpha _{1}\int_{\Omega \times Y_{1}}\left( \nabla p_{1}+\nabla _{y}\hat{p}%
_{1}\right) \mathbf{v}\left( x\right) dxdy+\alpha _{2}\int_{\Omega \times
Y_{2}}\nabla _{y}p_{2}\mathbf{v}\left( x\right) dxdy  \notag \\
&=&\alpha _{1}|Y_{1}|\int_{\Omega }\nabla p_{1}\mathbf{v}\left( x\right)
dx+\int_{\Omega \times \Gamma }\left( \alpha _{1}\hat{p}_{1}+\alpha
_{2}p_{2}\right) \left( \mathbf{v\cdot n}\right) dxds,  \label{37}
\end{eqnarray}%
By Theorem \ref{t1}, it follows that
\begin{eqnarray}
\lim_{\varepsilon \rightarrow 0}\int_{\Omega }\mathbf{f}^{\varepsilon }%
\mathbf{v}^{\varepsilon }\left( x\right) dx &=&\lim_{\varepsilon \rightarrow
0}\left( \int_{\Omega }\mathbf{f}^{\varepsilon }\left( x\right) \mathbf{v}%
\left( x\right) dx+\varepsilon \int_{\Omega }\mathbf{f}^{\varepsilon }\left(
x\right) \mathbf{\hat{v}}\left( x,\frac{x}{\varepsilon }\right) dx\right)
\notag \\
&=&\int_{\Omega }\mathbf{fv}\left( x\right) dx  \label{403}
\end{eqnarray}%
where $\mathbf{f}$ is given by (\ref{n6}). Thus, collecting these limits (%
\ref{401})-(\ref{403}), we obtain the limiting equation of (\ref{40})
\begin{eqnarray}
&&\int_{\Omega \times Y}\mathbb{A}\left[ \mathrm{e}\left( \mathbf{u}\right) +%
\mathrm{e}_{y}\left( \mathbf{\hat{u}}\right) \right] \left[ \mathrm{e}\left(
\mathbf{v}\right) +\mathrm{e}_{y}\left( \mathbf{\hat{v}}\right) \right]
dxdy+\alpha _{1}|Y_{1}|\int_{\Omega }\nabla p_{1}\mathbf{v}dx  \notag \\
&&+\int_{\Omega \times \Gamma }\left( \alpha _{1}\hat{p}_{1}+\alpha
_{2}p_{2}\right) \left( \mathbf{v\cdot n}\right) dxds=\int_{\Omega }\mathbf{%
fv}dx  \label{45}
\end{eqnarray}%
which is valid for a.e. $t\in \left( 0,T\right) $. Next, we proceed to get
the limiting equation of (\ref{24}). Taking $q_{1}=q_{1}^{\varepsilon }$ and
$q_{2}=q_{2}^{\varepsilon }$ in (\ref{24}), integrating by parts over $%
\left( 0,T\right) $ and taking into account the initial conditions (\ref{241}%
), we obtain
\begin{eqnarray}
&&-\int_{Q_{1}^{\varepsilon }}\left( c_{1}^{\varepsilon }\left( x\right)
p_{1}^{\varepsilon }+\alpha _{1}\mathrm{div}\mathbf{u}^{\varepsilon }\right)
\partial _{t}\varphi _{1}\left( t,x\right) dtdx-\int_{Q_{2}^{\varepsilon
}}c_{2}^{\varepsilon }\left( x\right) p_{2}^{\varepsilon }\partial
_{t}\varphi _{2}\left( t,x,\frac{x}{\varepsilon }\right) dtdx+  \notag \\
&&\int_{Q_{1}^{\varepsilon }}K_{1}\left( \frac{x}{\varepsilon }\right)
\nabla p_{1}^{\varepsilon }\left( \nabla \varphi _{1}\left( t,x\right)
+\nabla _{y}\hat{\varphi}_{1}\left( t,x,\frac{x}{\varepsilon }\right)
\right) dtdx+  \notag \\
&&\int_{Q_{2}^{\varepsilon }}\varepsilon k_{2}\left( \frac{x}{\varepsilon }%
\right) \nabla p_{2}^{\varepsilon }\nabla _{y}\varphi _{2}\left( t,x,\frac{x%
}{\varepsilon }\right) dtdx+  \notag \\
&&\varepsilon \int_{\Sigma ^{\varepsilon }}g\left( \frac{x}{\varepsilon }%
\right) \left( p_{1}^{\varepsilon }-p_{2}^{\varepsilon }\right) \left(
\varphi _{1}\left( t,x\right) -\varphi _{2}\left( t,x,\frac{x}{\varepsilon }%
\right) \right) dtds^{\varepsilon }+\varepsilon R_{2}^{\varepsilon }
\label{46}
\end{eqnarray}%
where
\begin{eqnarray*}
&&R_{2}^{\varepsilon }=\int_{Q_{1}^{\varepsilon }}-\left( c_{1}^{\varepsilon
}\left( x\right) p_{1}^{\varepsilon }+\alpha _{1}\mathrm{div}\mathbf{u}%
^{\varepsilon }\right) \partial _{t}\hat{\varphi}_{1}\left( t,x,\dfrac{x}{%
\varepsilon }\right) dtdx+ \\
&&\int_{Q_{2}^{\varepsilon }}-\alpha _{2}\mathrm{div}\mathbf{u}^{\varepsilon
}\partial _{t}\varphi _{2}\left( t,x,\frac{x}{\varepsilon }\right) dtdx+ \\
&&\int_{Q_{1}^{\varepsilon }}K_{1}\left( \frac{x}{\varepsilon }\right)
\nabla p_{1}^{\varepsilon }\nabla _{x}\hat{\varphi}_{1}\left( t,x,\frac{x}{%
\varepsilon }\right) dtdx+ \\
&&\varepsilon \int_{Q_{1}^{\varepsilon }}K_{2}\left( \frac{x}{\varepsilon }%
\right) \nabla p_{2}^{\varepsilon }\nabla _{x}\varphi _{2}\left( t,x,\frac{x%
}{\varepsilon }\right) dtdx+ \\
&&\varepsilon \int_{\Sigma ^{\varepsilon }}g\left( \frac{x}{\varepsilon }%
\right) \left( p_{1}^{\varepsilon }-p_{2}^{\varepsilon }\right) \hat{\varphi}%
_{1}\left( t,x\right) dtds^{\varepsilon }\text{.}
\end{eqnarray*}

The first integral of (\ref{46}) is equal to
\begin{equation*}
\int_{\Omega _{T}}-\chi _{1}\left( \frac{x}{\varepsilon }\right) \left(
c_{1}\left( \frac{x}{\varepsilon }\right) p_{1}^{\varepsilon }+\alpha _{1}%
\mathrm{div}\mathbf{u}^{\varepsilon }\right) \partial _{t}\varphi _{1}\left(
t,x\right) dtdx\text{,}
\end{equation*}%
and thanks to (\ref{32}) and (\ref{33}), converges to
\begin{equation*}
\int_{Q\times Y}-\chi _{1}\left( y\right) \left( c_{1}\left( y\right)
p_{1}+\alpha _{1}\left( \mathrm{div}\mathbf{u}+\mathrm{div}_{y}\mathbf{\hat{u%
}}\right) \right) \partial _{t}\varphi _{1}\left( t,x\right) dtdxdy.
\end{equation*}%
In a similar way, by (\ref{321}) and (\ref{33}), it follows that
\begin{equation*}
\int_{Q_{2}^{\varepsilon }}c_{2}^{\varepsilon }\left( x\right)
p_{2}^{\varepsilon }\partial _{t}\varphi _{2}\left( t,x,\frac{x}{\varepsilon
}\right) dtdx\rightarrow \int_{Q\times Y}\chi _{2}\left( y\right)
c_{2}\left( y\right) p_{2}\partial _{t}\varphi _{2}\left( t,x,y\right) dtdxdy
\end{equation*}%
Now, in view of (\ref{34}) one can deduce that
\begin{eqnarray*}
&&\int_{Q_{1}^{\varepsilon }}K_{1}\left( \frac{x}{\varepsilon }\right)
\nabla p_{1}^{\varepsilon }\left( \nabla \varphi _{1}\left( t,x\right)
+\nabla _{y}\hat{\varphi}_{1}\left( t,x,\frac{x}{\varepsilon }\right)
\right) dtdx= \\
&&\int_{Q}\chi _{1}\left( \frac{x}{\varepsilon }\right) K_{1}\left( \frac{x}{%
\varepsilon }\right) \nabla p_{1}^{\varepsilon }\left( \nabla \varphi
_{1}\left( t,x\right) +\nabla _{y}\hat{\varphi}_{1}\left( t,x,\frac{x}{%
\varepsilon }\right) \right) dtdx\rightarrow \\
&&\int_{Q\times Y}\chi _{1}\left( y\right) K_{1}\left( y\right) \left(
\nabla p_{1}+\nabla _{y}\hat{p}_{1}\right) \left( \nabla \varphi \left(
t,x\right) +\nabla _{y}\hat{\varphi}_{1}\left( t,x,y\right) \right) dtdxdy
\end{eqnarray*}%
and thanks to (\ref{35}), we also get
\begin{eqnarray*}
&&\int_{Q_{2}^{\varepsilon }}\varepsilon k_{2}\left( \frac{x}{\varepsilon }%
\right) \nabla p_{2}^{\varepsilon }\nabla _{y}\varphi _{2}\left( t,x,\frac{x%
}{\varepsilon }\right) dtdx \\
&=&\int_{Q}\chi _{2}\left( \frac{x}{\varepsilon }\right) K_{2}\left( \frac{x%
}{\varepsilon }\right) \varepsilon \nabla p_{2}^{\varepsilon }\nabla
_{y}\varphi _{2}\left( t,x,\frac{x}{\varepsilon }\right) dtdx\rightarrow \\
&&\int_{Q\times Y}\chi _{2}\left( y\right) K_{2}\left( y\right) \nabla
p_{2}\nabla _{y}\varphi _{2}\left( t,x,y\right) dtdxdy.
\end{eqnarray*}%
By virtue of (\ref{36}), we find that
\begin{eqnarray*}
&&\varepsilon \int_{\Sigma ^{\varepsilon }}g\left( \frac{x}{\varepsilon }%
\right) \left( p_{1}^{\varepsilon }-p_{2}^{\varepsilon }\right) \left(
\varphi _{1}\left( t,x\right) -\varphi _{2}\left( t,x,\frac{x}{\varepsilon }%
\right) \right) dtds^{\varepsilon }\rightarrow \\
&&\int_{Q\times \Gamma }g\left( y\right) \left( p_{1}-p_{2}\right) \left(
\varphi _{1}\left( t,x\right) -\varphi _{2}\left( t,x,y\right) \right)
dtdsdy.
\end{eqnarray*}

As before, we observe that $R_{2}^{\varepsilon }=O\left( 1\right) $ and, by
collecting all the preceeding limits, we get the following limiting equation
of (\ref{24}) :%
\begin{eqnarray}
&&\int_{Q\times Y_{1}}-\left( c_{1}\left( y\right) p_{1}+\alpha _{1}\left(
\mathrm{div}\mathbf{u}+\mathrm{div}_{y}\mathbf{\hat{u}}\right) \right)
\partial _{t}\varphi _{1}dtdxdy+  \notag \\
&&\int_{Q\times Y_{1}}K_{1}\left( y\right) \left( \nabla p_{1}+\nabla _{y}%
\hat{p}_{1}\right) \left( \nabla \varphi _{1}+\nabla _{y}\hat{\varphi}%
_{1}\right) dtdxdy+  \notag \\
&&\int_{Q\times Y_{2}}\left( -c_{2}\left( y\right) p_{2}\partial _{t}\varphi
_{2}+K_{2}\left( y\right) \nabla _{y}p_{2}\nabla _{y}\varphi _{2}\right)
dtdxdy+  \notag \\
&&\int_{Q\times \Gamma }g\left( y\right) \left( p_{1}-p_{2}\right) \left(
\varphi _{1}-\varphi _{2}\right) dtdsdy=0\text{.}  \label{50}
\end{eqnarray}

By density argument, the equations (\ref{45}) and (\ref{50}) still hold true
for any $\left( \mathbf{v},\mathbf{\hat{v}}\right) \in \mathbf{H}\times
L^{2}\left( \Omega ,H^{1}\left( Y\right) /\mathbb{R}\right) ^{3}$ and $%
\left( \varphi _{1},\hat{\varphi}_{1},\varphi _{2}\right) \in
L_{T}^{2}\left( H^{1}\left( \Omega \right) \right) \times L^{2}\left(
Q;H_{\#}^{1}\left( Y\right) /\mathbb{R}\right) \times L^{2}\left(
Q;H_{\#}^{1}\left( Y\right) \right) $. We can summarize the preceding by
observing that these equations are a weak formulation associated to the
two-scale homogenized system (\ref{471})-(\ref{484}). Indeed, integrating by
parts in (\ref{45}) and (\ref{50}), we obtain the following system:
\begin{eqnarray}
&&-\mathrm{div}_{y}\left( \mathbb{A}_{1}\left[ \mathrm{e}\left( \mathbf{u}%
\right) +\mathrm{e}_{y}\left( \mathbf{\hat{u}}\right) \right] \right) \text{%
$=$%
\ }0\text{ a.e.in }Q\times Y_{1},  \label{471} \\
&&\   \notag \\
&&-\mathrm{div}_{y}\left( \mathbb{A}_{2}\left[ \mathrm{e}\left( \mathbf{u}%
\right) +\mathrm{e}_{y}\left( \mathbf{\hat{u}}\right) \right] \right) \text{%
$=$%
\ }0\text{ a.e.in }Q\times Y_{2},  \label{472} \\
&&\   \notag \\
&&-\mathrm{div}\left( \int_{Y}\mathbb{A}\left[ \mathrm{e}\left( \mathbf{u}%
\right) +\mathrm{e}_{y}\left( \mathbf{\hat{u}}\right) \right] dy\right)
+\alpha _{1}\left\vert Y_{1}\right\vert \nabla p_{1}+  \notag \\
&&\int_{\Gamma }\left( \alpha _{1}\hat{p}_{1}+\alpha _{2}p_{2}\right)
\mathbf{n}ds\text{%
$=$%
\ }\mathbf{f}\text{ a.e. in }Q,  \label{473}
\end{eqnarray}%
and
\begin{eqnarray}
&&-\mathrm{div}_{y}\left( K_{1}\left( \nabla p_{1}+\nabla _{y}\hat{p}%
_{1}\right) \right) \text{%
$=$%
\ }0\text{ a.e. in }Q\times Y_{1},  \label{474} \\
&&\   \notag \\
&&\partial _{t}\left( c_{2}p_{2}\right) -\mathrm{div}_{y}\left( K_{2}\nabla
_{y}p_{2}\right) \text{%
$=$%
\ }0\text{ a.e.in }Q\times Y_{2},  \label{475} \\
&&\   \notag \\
&&\partial _{t}\left( \int_{Y_{1}}\left( c_{1}p_{1}+\alpha _{1}\left(
\mathrm{div}\mathbf{u}+\mathrm{div}_{y}\mathbf{\hat{u}}\right) \right)
\right) -\mathrm{div}\left( \int_{Y_{1}}K_{1}\left( \nabla p_{1}+\nabla _{y}%
\hat{p}_{1}\right) dy\right) +  \notag \\
&&\int_{\Gamma }g\left( y\right) \left[ p_{1}-p_{2}\right] ds\left( y\right)
\text{%
$=$%
\ }0\text{ a.e. in }Q,  \label{476}
\end{eqnarray}%
with the transmission and boundary conditions:%
\begin{eqnarray}
&&\mathbb{A}_{1}\left[ \mathrm{e}\left( \mathbf{u}\right) +\mathrm{e}%
_{y}\left( \mathbf{\hat{u}}\right) \right] \cdot \mathbf{n}\text{%
$=$%
\ }\mathbb{A}_{2}\left[ \mathrm{e}\left( \mathbf{u}\right) +\mathrm{e}%
_{y}\left( \mathbf{\hat{u}}\right) \right] \cdot \mathbf{n}\text{ a.e. on }%
Q\times \Gamma ,  \label{477} \\
&&\   \notag \\
&&\left( K_{1}\left( \nabla p_{1}+\nabla _{y}\hat{p}_{1}\right) \right)
\cdot \mathbf{n}\text{%
$=$%
\ }0\text{ a.e. on }Q\times \Gamma ,  \label{478} \\
&&\   \notag \\
&&  \notag \\
&&\ \ \left( K_{1}\left( \nabla p_{1}+\nabla _{y}\hat{p}_{1}\right) \right)
\cdot v\text{%
$=$%
\ }0\text{ a.e. on }\left( 0,T\right) \times \partial \Omega \times Y_{1},
\label{4781} \\
&&K_{2}\nabla _{y}p_{2}\cdot \mathbf{n}\text{%
$=$%
\ }-g\left( y\right) \left[ p_{1}-p_{2}\right] \text{ a.e. on }Q\times
\Gamma ,  \label{479} \\
&&\   \notag \\
&&\mathbf{u}\text{%
$=$%
\ }0\text{ a.e. on }\partial \Omega ,  \label{480} \\
&&\   \notag \\
&&y\longmapsto \ \mathbf{\hat{u}},\ \hat{p}_{1},p_{2}\text{ }Y\text{%
-periodic,}  \label{481}
\end{eqnarray}%
and the initial conditions:
\begin{eqnarray}
&&\mathbf{u}\left( 0,x\right) \text{%
$=$%
\ }\mathbf{0}\text{\ a.e. in }\Omega \text{,}  \label{482} \\
&&\mathbf{\hat{u}}\left( 0,x,y\right) \text{%
$=$%
\ }\mathbf{0}\text{\ a.e. in }\Omega \times Y\text{,}  \label{4821} \\
&&p_{1}\left( 0,x\right) \text{%
$=$%
\ }0\text{ a.e. in }\Omega \text{,}  \label{483} \\
&&\hat{p}_{1}\left( 0,x,y\right) \text{%
$=$%
\ }0\text{ a.e. in }\Omega \times Y_{1} \\
&&p_{2}\left( 0,x,y\right) \text{%
$=$%
\ }0\text{ a.e. in }\Omega \times Y_{2}\text{. }  \label{484}
\end{eqnarray}%
Now we decouple the system (\ref{471})-(\ref{484}). In view of the linearity
of the two first equations (\ref{471})-(\ref{472}), we can write that, up to
an additive constant:
\begin{equation}
\mathbf{\hat{u}}(t,x,y)=\sum_{i,j=1}^{3}\mathrm{e}_{ij}\left( \mathbf{u}%
\right) \left( t,x\right) \mathbf{w}^{ij}\left( y\right) +C^{te},\ \text{%
a.e. }\left( t,x,y\right) \in Q\times Y\text{, }  \label{48}
\end{equation}%
where, for $i,j\in \left\{ 1,2,3\right\} $,$\ \mathbf{w}^{ij}\in \left(
H_{\#}^{1}\left( Y\right) /\mathbb{R}\right) ^{3}$ is the solution to the
following microscopic system:
\begin{eqnarray*}
&&-\mathrm{div}_{y}\left( \mathbb{A}_{1}\mathrm{e}_{y}\left( \mathbf{w}^{ij}+%
\mathbf{d}^{ij}\right) \right) =0\text{ a.e. in }Y_{1},\ \ \ \  \\
&&\  \\
&&-\mathrm{div}_{y}\left( \mathbb{A}_{2}\mathrm{e}_{y}\left( \mathbf{w}^{ij}+%
\mathbf{d}^{ij}\right) \right) =0\text{ a.e. in }Y_{2}, \\
&&\  \\
&&\mathbb{A}_{1}\mathrm{e}_{y}\left( \mathbf{w}^{ij}+\mathbf{d}^{ij}\right)
\cdot \mathbf{n}=\mathbb{A}_{2}\mathrm{e}_{y}\left( \mathbf{w}^{ij}+\mathbf{d%
}^{ij}\right) \cdot \mathbf{n}\text{ a.e. on }\Gamma , \\
&& \\
&&\ y\longmapsto \ \mathbf{w}^{ij}\ Y\text{-periodic.}
\end{eqnarray*}%
Here $\mathbf{d}^{kl}=\left( y_{K}\delta _{il}\right) _{1\leq i\leq 3}$ and $%
\left( \delta _{ij}\right) $ is the Kr\"{o}necker symbol.

Similarly, in view of (\ref{474}), (\ref{478}) and (\ref{481}) one can write
that:
\begin{equation}
\hat{p}_{1}\left( t,x,y\right) =\sum_{i=1}^{3}\frac{\partial p_{1}}{\partial
x_{i}}\left( t,x\right) \pi _{i}\left( y\right) +C^{te},\ \text{a.e. }\left(
t,x,y\right) \in Q\times Y_{1}\text{,}  \label{49}
\end{equation}%
where, for $i=1,2,3$, the micro-pressure $\pi _{i}\in H^{1}\left(
Y_{1}\right) /\mathbb{R}$ is the solution of the following stationary
equation:
\begin{gather*}
-\mathrm{div}_{y}\left( K_{1}\left( \nabla \pi _{i}+e_{i}\right) \right) =0%
\text{ in }Y_{1}\text{,} \\
K_{1}\left( \nabla \pi _{i}+e_{i}\right) \cdot \mathbf{n}=0\text{ on }\Gamma
\text{,} \\
y\longmapsto \pi _{i}\text{ }Y\text{-periodic.}
\end{gather*}%
Here $e_{i}$ is the $i^{\text{th}}$ vector of the canonical basis of $%
\mathbb{R}^{3}$. Let us denote
\begin{eqnarray*}
\widetilde{\mathbb{A}} &=&\left( \tilde{a}_{i_{1}i_{2}i_{3}i_{4}}\right)
_{1\leq i_{1},i_{2},i_{3},i_{4}\leq 3},\ \  \\
\tilde{a}_{i_{1}i_{2}i_{3}i_{4}}
&=&\sum_{j_{1},j_{2}=1}^{3}\int_{Y}a_{i_{1}i_{2}j_{1}j_{2}}\left( y\right)
\left( \delta _{i_{1}j_{1}}\delta _{i_{2}j_{2}}+\mathrm{e}%
_{j_{1}j_{2},y}\left( \mathbf{w}^{i_{3}i_{4}}\right) \left( y\right) \right)
dy,\ \
\end{eqnarray*}%
where $\left( a_{ijlm}\right) $ are the coefficients of the elasticity
tensor $\mathbb{A}$ and
\begin{equation*}
\mathrm{e}_{ij,y}\left( \mathbf{w}\right) =\frac{1}{2}\left( \frac{\partial
w_{i}}{\partial y_{j}}+\frac{\partial w_{j}}{\partial y_{i}}\right) ,\
\mathbf{w}=\left( w_{j}\right) _{1\leq j\leq 3}\text{.}
\end{equation*}%
Let also define the effective stress tensor
\begin{equation*}
\tilde{\sigma}(\mathbf{u})=\left( \tilde{\sigma}_{ij}(\mathbf{u})\right)
_{1\leq i,j\leq 3},\ \ \tilde{\sigma}_{ij}(\mathbf{u})=\sum_{l,m=1}^{3}%
\tilde{a}_{ijlm}\mathrm{e}_{lm}(\mathbf{u})\text{,\ }
\end{equation*}%
the effective permeability tensor
\begin{equation*}
\tilde{K}=\left( \tilde{K}_{ij}\right) _{1\leq i,j\leq 3},\ \ \tilde{K}%
_{ij}=\int_{Y_{1}}K_{1}\left( y\right) \left( \nabla _{y}\pi
_{i}+e_{i}\right) \left( \nabla \pi _{j}+e_{j}\right) dy,
\end{equation*}%
the effective Biot-Willis matrices:%
\begin{eqnarray*}
B &=&\left( b_{ij}\right) ,\ b_{ij}=\alpha _{1}\left( \left\vert
Y_{1}\right\vert \delta _{ij}+\int_{\Gamma }\pi _{j}\left( y\right)
n_{i}ds\left( y\right) \right) ,\ \mathbf{n}=\left( n_{i}\right) _{1\leq
i\leq 3} \\
\Lambda &=&\left( \lambda _{ij}\right) _{1\leq i,j\leq 3},\ \ \lambda
_{ij}=\alpha _{1}\int_{Y_{1}}\sum_{m=1}^{3}\left( \delta _{im}\delta _{jm}+%
\frac{\partial w_{m}^{ij}}{\partial y_{m}}\right) dy, \\
\mathbf{w}^{ij} &=&\left( w_{m}^{ij}\right) _{1\leq m\leq 3}
\end{eqnarray*}%
and finally the following averaging quantities%
\begin{equation*}
\tilde{c}=\int_{Y_{1}}c_{1}\left( y\right) dy,\ \ \ \tilde{g}=\int_{\Gamma
}g\left( y\right) ds\left( y\right) \text{.}
\end{equation*}%
Then from (\ref{48})-(\ref{49}) we deduce the homogenized system :
\begin{gather}
-\mathrm{div}\tilde{\sigma}(\mathbf{u})+B\nabla p_{1}+\alpha
_{2}\int_{\Gamma }p_{2}\mathbf{n}ds\left( y\right) \text{%
$=$%
\ }\mathbf{f}\text{ a.e. in }Q\text{,}  \label{d1} \\
\partial _{t}\left( \tilde{c}p_{1}+\Lambda :\mathrm{e}\left( \mathbf{u}%
\right) \right) -\mathrm{div}\left( \tilde{K}\nabla p_{1}\right) +\tilde{g}%
p_{1}  \notag \\
-\int_{\Gamma }g\left( y\right) p_{2}ds\left( y\right) =0,\text{ a.e. in }Q,
\label{d2} \\
\partial _{t}\left( c_{2}p_{2}\right) -\mathrm{div}_{y}\left( K_{2}\nabla
_{y}p_{2}\right) \text{%
$=$%
\ }0\text{ a.e. in }Q\times Y_{2},  \label{d3} \\
c_{2}\nabla _{y}p_{2}\cdot \mathbf{n}\text{%
$=$%
\ }-g\left( y\right) \left[ p_{1}-p_{2}\right] \text{ a.e. on }Q\times
\Gamma ,  \label{d31} \\
\mathbf{u}\text{%
$=$%
\ }0\text{ and }\tilde{K}\nabla p_{1}\cdot \nu =0\text{ a.e. on }\left(
0,T\right) \times \Sigma ,  \label{d4} \\
y\longmapsto \ p_{2}\text{ }Y\text{-periodic,}  \label{d5} \\
\mathbf{u}\left( 0,x\right) \text{%
$=$%
\ }\mathbf{0}\text{\ a.e. in }\Omega \text{, }p_{1}\left( 0,x\right) \text{%
$=$%
\ }0\text{ a.e. in }\Omega \text{,}  \label{d6} \\
p_{2}\left( 0,x,y\right) \text{%
$=$%
\ }0\text{ a.e. in }\Omega \times Y_{2}\text{.}  \label{d9}
\end{gather}

Now, we are going to establish a relation between the two pressures $p_{1}$
and $p_{2}$. To this aim, let $\zeta \in L^{\infty }\left(
0,T;H_{\#}^{1}\left( Y_{2}\right) \right) $ be the unique solution to the
following microscopic and non homogeneous Robin problem:%
\begin{eqnarray*}
&&\partial _{t}\left( c_{2}\zeta \right) -\mathrm{div}_{y}\left( K_{2}\nabla
_{y}\zeta \right) \text{%
$=$%
\ }0\text{ a.e.in }\left( 0,T\right) \times Y_{2}\text{,} \\
&&K_{2}\nabla _{y}\zeta \cdot \mathbf{n}\text{%
$=$%
\ }-g\left( y\right) \left[ 1-\zeta \right] \text{ a.e. on }\Sigma \text{,}
\\
&&y\longmapsto \zeta \text{ }Y\text{-periodic,} \\
&&\zeta \left( 0,y\right) =0\text{,\ a.e. }y\in Y_{2}\text{.}
\end{eqnarray*}%
Since $c_{2},K_{2},g$ are time-independent and $p_{1}$ is independent of $y$%
, using the Laplace transform method, one can then easily see that
\begin{equation}
p_{2}\left( t,x,y\right) =\int_{0}^{t}p_{1}\left( \tau ,x\right) \partial
_{t}\zeta \left( t-\tau ,y\right) d\tau ,\ \text{a.e. }\left( t,x,y\right)
\in Q\times Y_{2}\text{.}  \label{rel1}
\end{equation}%
Therefore, the homogenized system (\ref{d1})-(\ref{d9}) can be rewritten as
\begin{eqnarray*}
&&-\mathrm{div}\tilde{\sigma}(\mathbf{u})+B\nabla p_{1}+\int_{0}^{t}\theta
\left( t,\tau \right) p_{1}\left( \tau ,x\right) d\tau \text{%
$=$%
\ }\mathbf{f}\text{ a.e. in }Q\text{,} \\
&&\partial _{t}\left( \tilde{c}p_{1}+\Lambda :\mathrm{e}\left( \mathbf{u}%
\right) \right) -\mathrm{div}\left( \tilde{K}\nabla p_{1}\right) +\tilde{g}%
p_{1} \\
&&-\int_{0}^{t}\eta \left( t,\tau \right) p_{1}\left( \tau ,x\right) d\tau
\text{%
$=$%
}0,\text{ a.e. in }Q, \\
&&\mathbf{u}\text{%
$=$%
\ }0\text{ and }\tilde{K}\nabla p_{1}\cdot \nu \text{%
$=$
}0\text{ a.e. on }\left( 0,T\right) \times \partial \Omega , \\
&&\mathbf{u}\left( 0,x\right) \text{%
$=$%
\ }\mathbf{0}\text{,\ }p_{1}\left( 0,x\right) \text{%
$=$%
\ }0\text{ a.e. in }\Omega \text{,}
\end{eqnarray*}%
where we have denoted
\begin{eqnarray*}
\theta \left( t,\tau \right) &=&\alpha _{2}\int_{\Gamma }\partial _{t}\zeta
\left( t-\tau ,y\right) \mathbf{n}ds\left( y\right) ,\  \\
\eta \left( t,\tau \right) &=&\int_{\Gamma }g\left( y\right) \partial
_{t}\zeta \left( t-\tau ,y\right) ds\left( y\right) \text{.}
\end{eqnarray*}%
Finally, let us observe that the overall pressure of the fluid flow in the
microstructure model which is given by
\begin{equation*}
P^{\varepsilon }\left( t,x\right) =\chi _{1}^{\varepsilon
}(x)p_{1}^{\varepsilon }\left( t,x\right) +\chi _{2}^{\varepsilon
}(x)p_{2}^{\varepsilon }\left( t,x\right)
\end{equation*}%
for a.e. $\left( t,x\right) \in Q$, two-scale converges to $\chi
_{1}(y)p_{1}\left( t,x\right) +\chi _{2}(y)p_{2}\left( t,x,y\right) $, and
thanks to (\ref{rel1}), converges then weakly in $L^{2}\left( Q\right) $ to%
\begin{equation*}
|Y_{1}|p_{1}\left( t,x\right) +\int_{0}^{t}\int_{Y_{2}}p_{1}\left( \tau
,x\right) \partial _{t}\zeta \left( t-\tau ,y\right) dyd\tau \text{.}
\end{equation*}%
This concludes the proof of Theorem \ref{thp}.

\section{Conclusion\label{s4}}

We have used the homogenization theory to derive a macro-model for fluid
flow in composite poroelastic with microstructures, in which inclusions are
fully embedded and with very low permeabilities. We have shown that the
overall behavior of fluid flow in such heterogeneous media with low
permeability at the micro-scale may present memory terms. We also have shown
that in such cases, the Biot-Willis parameters are, as in \cite{ain2},
matrices and no longer scalars, as it is usually considered in the
poroelasticity litterature, since it is assumed there that the medium is
homogeneous and isotropic. Nevertheless, anisotropic media may present
different coupling interaction properties in different directions at the
micro-scale, and which lead at the macro-scale to such anisotropic
Biot-Willis parameters. Finally, let us mention that the result of the paper
remains valid if one considers non homogeneous initial conditions or with
any volume distributed source densities in each phases.

\end{document}